\documentclass[10pt,leqno,twoside]{article}

\usepackage[latin2]{inputenc}
\usepackage{amsmath}
\usepackage{amsfonts}
\usepackage{amssymb}
\usepackage{makeidx}
\newcommand{\RR}{{\mathbb R}}

\newcommand{\NN}{{\mathbb N}}

\newtheorem{statement}{Statement}
\newtheorem{theorem}{Theorem}

\newtheorem{lemma}{Lemma}

\newtheorem{remark}[theorem]{Remark}

\newfont{\eurorm}{eurm10 scaled 1100}
\newfont{\eurorms}{eurm10 scaled 800}

\date{April 8, 2018}

\begin{document}

\bigskip
\title{\bf About some exponential inequalities \\ related to the sinc function}
\maketitle

\begin{center}
{\em Marija Ra\v sajski, Tatjana Lutovac, Branko Male\v sevi\' c${}^{\;\mbox{\scriptsize $\ast$}}\!$}
\end{center}

\begin{center}
{\footnotesize
\textit{
School of Electrical Engineering, University of Belgrade, \\[-0.25 ex]
Bulevar kralja Aleksandra 73, 11000 Belgrade, Serbia}}
\end{center}

\medskip
\noindent {\small \textbf{Abstract.}
{\small
In this paper we prove some exponential  inequalities involving the sinc function. We analyze and prove inequalities
with constant exponents as well as inequalities with certain polynomial exponents.$\;$Also, we establish intervals
in which these inequalities hold.}

\footnote{$\!\!\!\!\!\!\!\!\!\!\!\!\!\!$
{\scriptsize
${}^{\mbox{\scriptsize $\ast$}}$Corresponding author.       \\[0.0 ex]
Emails:                              \\[0.0 ex]
{\em Marija Ra\v sajski} {\tt $<$marija.rasajski@etf.rs$>$},
{\em Tatjana Lutovac} {\tt $<$tatjana.lutovac@etf.rs$>$},     \\[0.0 ex]
{\em Branko Male\v sevi\' c} {\tt $<$branko.malesevic@etf.rs$>$}}}

{\footnotesize Keywords: Exponential inequalities, sinc function }

{\small \tt MSC: Primary 33B10; Secondary 26D05}

\section{Introduction}

Inequalities related  to the  sinc function, i.e. $\displaystyle \mbox{\rm sinc}\,x\!=\!\mbox{\small $\dfrac{\sin x}{x}$}$
\mbox{${\big (}\displaystyle  x \neq 0 {\big )}$},
occur  in many fields of mathematics and  engineering \cite{D_S_Mitrinovic_1970}, \cite{Mortici_2011}, \cite{Rahmatollahi_DeAbreu_2012}, \cite{G_V_Milovanovic_2014}, \cite{Cloud_Drachman_Lebedev_2014}, \cite{RIM2018}, \cite{AIDE2018}  such as {\sc Fourier}\ analysis and its applications, information theory, radio transmission, optics, signal processing, sound recording, etc.


\medskip \noindent  The following
inequalities are proved in \cite{Z.-H._Yang_2014}:
\begin{equation}
\cos ^{2}{\!\displaystyle\frac{x}{2}}\leq \displaystyle\frac{\sin {x}}{x}\leq \cos ^{3}{\!\displaystyle\frac{x}{3%
}}\leq \displaystyle\frac{2+\cos {x}}{3} \label{Z-H-Jang-1}
\end{equation}%
for every $x\in \left( 0,\pi \right). $

\smallskip

\noindent In \cite{Lutovac2017}, the authors considered
 possible refinements of the inequality (\ref%
{Z-H-Jang-1}) by a real analytic function $\varphi
_{a}(x)\!=\!\left( \displaystyle\frac{\sin x}{x}\right) ^{\!a}\!\!,
$ for $x\!\in
\!\left( 0,\pi \right) $ and parameter $a\!\in \!%
\mathbb{R}
$, and proposed and proved the following inequalities:

\begin{statement}{\rm (\cite{Lutovac2017}, Theorem 10)}
\label{Brankova-teorema}
The following inequalities hold true, for every ${x \!\in\! \left(0, \pi\right)}$
and $a \!\in\! \displaystyle\left(1, \mbox{\footnotesize $\displaystyle\frac{3}{2}$}\right):$

\vspace*{-3.0 mm}

\begin{equation} \label{Z-H-Jang}
\cos^{2}{\!\displaystyle\frac{x}{2}} \leq
\left(\displaystyle\frac{\sin{x}}{x}\right)^{\!a} \leq
\displaystyle\frac{\sin{x}}{x}.
\end{equation}
\end{statement}

In the paper \cite{Lutovac2017}, based on the analysis of the sign of the analytic function
$$
F_a(x)
=
\left( \displaystyle\frac{\sin x}{x}\right) ^{\!a}-\cos ^{2}{\!\displaystyle\frac{x}{2}}
$$
in the right neighborhood of zero, the corresponding inequalities for parameter values
$\displaystyle a \geq \mbox{\footnotesize $\displaystyle\frac{3}{2}$}$ are discussed.

\medskip
In this paper, in subsection 3.1,  using the power series expansions and the application of the {\sc Wu}-{\sc Debnath} theorem, we prove that the inequality (\ref{Z-H-Jang})
holds for $a = \mbox{\footnotesize $\displaystyle\frac{3}{2}$}$. At the same time, this proof represents  another proof of  Statement~\ref{Brankova-teorema}.
Also, we analyze the cases $a\in \left(\mbox{\footnotesize $\displaystyle\frac{3}{2}$},2\right)$ and $\displaystyle a \geq 2$ and we prove  the corresponding inequalities.

In subsection 3.2 we introduce and prove a new double-sided inequality of similar type involving polynomial
exponents.

\smallskip

Finally, in subsection 3.3, we establish a relation between the cases of the constant  and of the polynomial  exponent.

\smallskip

\section{Preliminaries}

In this section we review some results that we use in our study.

\smallskip
In accordance with \cite{Gradshteyn-Ryzhik}, the following expansions hold:
\begin{equation}
\label{Series_ln_sin_x_over_x}
\ln \frac{\sin x}{x}
=
-\sum\limits_{k=1}^{\infty}{\frac{2^{2k-1}|\mbox{\bf B}_{2k}|}{k(2k)!}x^{2k}}, \qquad (0 < x < \pi),
\end{equation}
\begin{equation}
\label{Series_ln_cos_x}
\ln \cos x
=
-\sum\limits_{k=1}^{\infty}{\frac{2^{2k-1}(2^{2k}-1)|\mbox{\bf B}_{2k}|}{k(2k)!}x^{2k}}, \qquad (-\pi/2 < x < \pi/2),
\end{equation}
where $\mbox{\bf B}_{i}$ ($i \!\in\! \NN$) are {\sc Bernoulli}'s numbers.

\medskip
The following theorem proved by {\sc Wu} and {\sc Debnath} in \cite{Wu_Debnath_2009}, is used in our proofs.

\medskip
\noindent {\bf Theorem WD.} (\cite{Wu_Debnath_2009}, Theorem 2 )
\label{Debnath_Wu_T}
{\em Suppose that $f(x)$ is a real function on $(a,b)$, and that $n$ is a positive integer such that $f ^{(k)}(a+), f^{(k)}(b-)$,
$\left(k \!\in\! \{0,1,2, \ldots ,n\}\right)$ exist.

\medskip
\noindent
{\boldmath $(i)$} Supposing that $(-1)^{(n)} f^{(n)}(x)$~is~in\-cre\-asing on $(a,b)$, then
for all $x \in (a,b)$ the following inequality holds$\,:$
\begin{equation}
\label{Debnath_Wu_first}
\begin{array}{c}
\displaystyle\sum_{k=0}^{n-1}{\mbox{\small $\displaystyle\frac{f^{(k)}(b\mbox{\footnotesize $-$})}{k!}$}(x\!-\!b)^k}
+
\frac{1}{(a-b)^n}
{\bigg (}\!
f(a\mbox{\footnotesize $+$})
-
\displaystyle\sum_{k=0}^{n-1}{\mbox{\small $\displaystyle\frac{(a\!-\!b)^{k}f^{(k)}(b\mbox{\footnotesize $-$})}{k!}$}
\!{\bigg )} (x\!-\!b)^{n}}                                          \\[2.5 ex]
<
f(x)
<
\displaystyle\sum_{k=0}^{n}{\frac{f^{(k)}(b\mbox{\footnotesize $-$})}{k!}(x\!-\!b)^{k}}.
\end{array}
\end{equation}
Furthermore, if $(-1)^{n} f^{(n)}(x)$ is decreasing on $(a,b)$, then the reversed inequality of {\rm (\ref{Debnath_Wu_first})} holds.

\break

\medskip
\noindent
{\boldmath $(ii)$} Supposing that $f^{(n)}(x)$ is increasing on $(a,b)$, then for all $x \!\in\! (a,b)$
the following inequality also holds$\,:$
\begin{equation}
\label{Debnath_Wu_second}
\begin{array}{l}
\displaystyle\sum_{k=0}^{n}{\frac{f^{(k)}(a\mbox{\footnotesize $+$})}{k!}(x-a)^{k}}
 <
f(x) < \\[2.5 ex]
<
\displaystyle \sum_{k=0}^{n-1}{\mbox{\small $\displaystyle\frac{f^{(k)}(a\mbox{\footnotesize $+$})}{k!}$}(x\!-\!a)^k}
+
\frac{1}{(b\!-\!a)^n}
{\bigg (}\!
f(b-)
-
\displaystyle\sum_{k=0}^{n-1}{\mbox{\small $\displaystyle\frac{(b-a)^{k}f^{(k)}(a\mbox{\footnotesize $+$})}{k!}$}
\!{\bigg )} (x\!-\!a)^{n}}.
\end{array}
\end{equation}
Furthermore, if $f^{(n)}(x)$ is decreasing on $(a,b)$, then the reversed inequality~of~\mbox{\rm (\ref{Debnath_Wu_second})} holds.}
%
%
\begin{remark} Note that inequalities $(\ref{Debnath_Wu_first})$ and $ (\ref{Debnath_Wu_second})$ hold for $n \in \NN$ as well as for $n=0$. \\
Here, and throughout this paper, a sum where the upper bound of the summation is lower than the lower bound of the summation, is understood to be zero.
\end{remark}

\medskip
The following Theorem, which is a consequence of Theorem WD, was proved~in~\cite{JNSA2018}.

%
%

\begin{theorem}{\rm (\cite{JNSA2018}, Theorem 1)}
\label{Natural_Extension_Theorem}
Let the function $f\!:\!(a,b) \longrightarrow \RR$ have the following power series expansion$\,:$
\begin{equation}
f(x)
=
\displaystyle\sum_{k=0}^{\infty}{c_{k}(x-a)^k}
\end{equation}
for $x \!\in\! (a,b)$, where the sequence of coefficients $\{c_{k}\}_{k \in \NN_0}$  has
a finite number of non-positive members and their indices are in the set $J \!=\! \{j_0,\ldots,j_\ell\}$.

\smallskip
Then, for the function
\begin{equation}
F(x)
=
f(x)-\displaystyle\sum_{i=0}^{\ell}{c_{j_i}(x-a)^{j_i}}
=
\displaystyle\sum_{k \in \NN_0 \backslash\!\;\! J}{c_{k}(x-a)^k},
\end{equation}
and the sequence $\{C_{k}\}_{k \in \NN_0}$ of the non-negative coefficients defined by$:$
\begin{equation}
C_{k}
=
\left\{
\begin{array}{ccc}
c_{k} \!&\!:\!&\! c_{k} > 0,        \\[0.75 ex]
0   \!&\!:\!&\! c_{k} \leq 0;
\end{array}
\right.
\end{equation}
holds that$\,:$
\begin{equation}
F(x)
=
\displaystyle\sum_{k=0}^{\infty}{C_{k}(x-a)^k},
\end{equation}
for every $x \!\in\! (a,b)$.

\bigskip
It is also $F^{(k)}(a+)= k!\,C_{k}$ and the following inequalities hold$:$
\begin{equation}
\!
\begin{array}{l}
\displaystyle\sum_{k=0}^{n}{C_k(x-a)^{k}} < F(x) <                              \\[3.0 ex]
<
\displaystyle\sum_{k=0}^{n-1}{C_k(x-a)^k}
+
\frac{1}{(b-a)^n}
{\bigg (}
F(b-)
-
\displaystyle\sum_{k=0}^{n-1}{C_k(b-a)^k}
{\bigg )} (x-a)^{n},
\end{array}
\end{equation}
for every $x \in (a, b)$ and $n \in \NN_0$, i.e.

\begin{equation}
\!\!
\begin{array}{l}
\displaystyle \displaystyle\sum_{k=0}^{m}{\!C_k}{(x\!-\!a)^k}
+
\displaystyle\sum_{i=0}^{\ell}{\!c_{j_i}}{(x - a)^{j_i}} \, < \, f(x) \, <                  \\[3.0 ex]
< \,
\displaystyle
\sum_{k = 0}^{m - 1}{\!C_k}{(x \!-\! a)^{k}}
\!+\!
\displaystyle\sum_{i = 0}^{\ell}{\!c_{j_i}}{(x \!-\! a)^{j_i}}
\!+\!
\displaystyle\frac{{(x \!-\! a)}^m}{{(b \!-\! a)}^m}
\!\left(\!{f(b\mbox{\footnotesize $-$})
\!-\!
\displaystyle\sum_{k=0}^{m-1}{\!C_k}{{(b\!-\!a)}^k}
\!-\!
\displaystyle\sum_{i=0}^{\ell}{\!c_{j_i}}{(b\!-\!a)}^{j_i}}\!\right)
\end{array}
\end{equation}
for every $x \!\in\! (a,b)$ and $m > max\{j_0,\ldots,j_\ell\}$.
\end{theorem}

\smallskip

\section{Main results}
\subsection{Inequalities with  constants in the exponents}
First, we consider a connection between the number of zeros of a real analytic function and some properties of its derivatives.
It is well known that the zeros of a non-constant analytic function are isolated
\cite{Godement_2004}, see also \cite{Krantz_Parks_1992} and \cite{Malesevic2016}.

\medskip
We prove the following assertion:
\begin{theorem}
\label{exactly_one_zero}
Let $f\!: (0, c) \longrightarrow \RR$ be real analytic function  such that
$f^{(k)}(x) > 0$ for $x \in (0,c)$ and $k= m, m+1, \ldots\, , {\big (}\,$for some $m \in \NN\,{\big )}$.

\smallskip
\noindent\enskip
If the following conditions hold$:$
\begin{itemize}
\item[$1)$]
there is a right neighbourhood  of zero in which  the following inequalities hold true:
$\, f(x)<0, \, f'(x)<0, \, \ldots, f^{(m-1)}(x)<0,$
\item[and]
\item[$2)$] $f(c_-)>0, f'(c_-)>0, \ldots, f^{(m-1)}(c_-)>0,$
\end{itemize}
then there exists exactly one zero $ x_0 \in (0, c) $ of the function $ f $.

\end{theorem}
{\bf Proof.} As $f^{(m)}(x) \!>\! 0$ for $x \!\in\! (0,c)$, it follows that $f^{(m-1)}(x)$ is monotonically increasing function for
$x \!\in\! (0,c)$. Based on  conditions $1)$ and $2)$, we conclude that there exists exactly one zero $ x_{m-1} \!\in\! (0, c)$
of the function $f^{(m-1)}(x)$. Next, we can conclude that function $f^{(m-2)}(x)$ is monotonically decreasing for $x \!\in\! (0,x_{m-1})$
and monotonically increasing for $x \!\in\! (x_{m-1},c)$. It is clear that function $ f^{(m-2)}(x) $ has exactly one minimum in the interval
$(0, c)$ at point $x_{m-1}$ and $f^{(m-2)}(x_{m-1}) \!<\! 0$.  On the basis of conditions 2), it follows that function $f^{(m-2)}(x)$ has exactly
one root $x_{m-2}$ on the interval $(0, c)$  and $x_{m-2} \!\in\! (x_{m-1}, c)$.

\smallskip

\noindent By repeating the described procedure, we get the assertion given in the theorem.\hfill $\Box$

\bigskip
Let us consider the family of functions
\begin{equation}\displaystyle
f_a(x) = a \ln \frac{\sin x}{x} - 2 \ln \cos\frac{x}{2},
\end{equation}
for $x \!\in\! (0, \pi)$ and parameter $a \in (1, +\infty)$.


\smallskip
Obviously,  the following equivalence  is true:
\begin{equation}\label{f_a}
a_1 < a  \Longleftrightarrow f_{a}(x) < f_{a_1}(x),
\end{equation}
for $a, a_1  > 1 $ and  $x \!\in\! (0,\pi)$.

\break

Thus:
\begin{equation}
\frac{3}{2} < a \Longleftrightarrow f_{a}(x) < f_{\frac{3}{2}}(x), \quad \quad \mbox{for $x \!\in\! (0,\pi)$}.
\end{equation}

Based on the power series expansions (\ref{Series_ln_sin_x_over_x}) and (\ref{Series_ln_cos_x}), we have:
\begin{equation}
\label{fun_f}
f_a(x)
=
\displaystyle\sum_{k=1}^{\infty}{E_k \, x^{2k}}
\end{equation}
for $a>1$ and $x \in (0, \pi)$, where
\begin{equation}\label{power-series-f}
E_k = \displaystyle\frac{{\big (}(2-a)\,4^k-2\,{\big )}|\mbox{\bf B}_{2k}|}{2 k \cdot (2k)!} \quad (k \!\in\! \NN).
\end{equation}

For $a \!=\! \mbox{\footnotesize $\displaystyle\frac{3}{2}$}$,  it is true that $E_{1}=0$ and $E_{k}>0$  for $k=2,3 , \ldots$. Thus, from (\ref{fun_f}), we have
$$ \displaystyle f_{\frac{3}{2}}(x)> 0 \quad \quad \mbox{for} \,\, x \in (0, \pi),$$
and  consequently the following theorem holds:

%
%

\begin{theorem}  The following inequalities hold true, for every $x \in (0, \pi)\!:$
$$
\cos^{2}{\!\displaystyle\frac{x}{2}} \leq
\left(\displaystyle\frac{\sin{x}}{x}\right)^{\frac{3}{2}} \leq
\displaystyle\frac{\sin{x}}{x}.
$$
\end{theorem}

\smallskip
As the  inequality
$$
\left(\displaystyle\frac{\sin{x}}{x}\right)^{\frac{3}{2}} \leq \left(\displaystyle\frac{\sin{x}}{x}\right)^{a}
$$
holds for $x \!\in\! (0, \pi)$ and $a \!\in\! \left(1, \mbox{\footnotesize $\displaystyle\frac{3}{2}$}\right]$,
the previous theorem can be thought of as a new
proof of Statement~\ref{Brankova-teorema}.

\medskip

Consider now the family of functions
$
\displaystyle
f_a(x) = a \ln \frac{\sin x}{x} - 2 \ln \cos\frac{x}{2},
$
for $x \!\in\! (0, \pi)$ and parameter $a \!>\! \mbox{\footnotesize $\displaystyle\frac{3}{2}$}$.

\smallskip

It easy to check that for  the sequence $\{\alpha_k\}_{k \in \NN}$ where
\begin{equation}
\alpha_k=2-\frac{2}{4^k}.
\end{equation}
the following equivalences are true:
\begin{equation}
\begin{array}{l}
a = \alpha_k  \;\Longleftrightarrow\; E_k = 0, \\[0.75 em]
a \in \left(\alpha_k, \alpha_{k+1} \right)
\;\Longleftrightarrow\;
\left(\forall i \!\in\! \{1,2,\ldots,k \} \right) E_i < 0
\;\wedge\;
\left(\forall i \! > \! k \right) E_i > 0.
\end{array}
\end{equation}

Let us now consider the function $\mathfrak{m} \!:\! \left[\mbox{\footnotesize $\displaystyle\frac{3}{2}$},2\right) \longrightarrow \NN_0$ defined by:
\begin{equation}\label{definition_m(a)}
\mathfrak{m}(a)=k \;\;\;\mbox{if and only if }\;\; a \in \left(\alpha_k, \alpha_{k+1} \right].
\end{equation}

\smallskip

It is not difficult to check  that $\lim\limits_{a \rightarrow 2_{-}}\!\!\mathfrak{m}(a) = +\infty$,
while for a fixed $a \!\in\! \left(\mbox{\footnotesize $\displaystyle\frac{3}{2}$},2\right)$ the number of negative elements of the sequence $\{ E_k \}_{k \in \NN}$ is $\mathfrak{m}(a)$ and their indices are in the set $\{ 1, \ldots , \mathfrak{m}(a)\}$.  For this reason, we distinguish two cases  $a \!\in\! \left(\mbox{\footnotesize $\displaystyle\frac{3}{2}$},2\right)$ or $\,a\geq 2$.

 As for the parameter $a = 2$ and $x \in (0, \pi)$  we have:
$$  \left(\displaystyle\frac{\sin{x}}{x}\right)^{\!2} \leq \cos^{2}{\!\displaystyle\frac{x}{2}}   \, \Longleftrightarrow \, \displaystyle \sin^{2}\dfrac{x}{2} \leq  \left(\dfrac{x}{2}\right)^{2},$$
while
for $a>2$ and $x \in (0, \pi)$  we have:
$$ \left(\displaystyle\frac{\sin{x}}{x}\right)^{\!a} \leq \left(\displaystyle\frac{\sin{x}}{x}\right)^{\!2}.$$
Hence, we have proved the following theorem:

\begin{theorem}
For every $a \geq 2$ and every $x \in  (0, \pi)$
the following inequality  holds true$:$
\begin{equation}
\label{nejednakost_Corollary11}
\left(\displaystyle\frac{\sin{x}}{x}\right)^{\!a} \leq \cos^{2}{\!\displaystyle\frac{x}{2}}.
\end{equation}
\end{theorem}

\medskip
 Consider now the  case  when the parameter $a \in \left(\mbox{\footnotesize $\displaystyle\frac{3}{2}$},2\right) $.
As noted above, for any fixed $a \in \left(\mbox{\footnotesize $\displaystyle\frac{3}{2}$},2\right) $ there is a finite number of negative
coefficients in the power series expansion~(\ref{power-series-f}), so it is possible to apply Theorem~\ref{Natural_Extension_Theorem}.

\smallskip

According to Theorem~\ref{Natural_Extension_Theorem}, the following inequalities  hold:

\vspace*{-2.5 mm}

%
\begin{equation}\label{natural-extension-estimation}
\begin{array}{l}
\displaystyle\sum_{\,\,\,k=\mathfrak{m}(a)+1}^{n}{\!\!\!E_k}{x^k}
+
\!\!
\displaystyle\sum_{i=0}^{\mathfrak{m}(a)-1}{\!E_{i}}{x^{i}} \, <                       \\[1.75 em]
< \, f_a(x) \, <                                                                       \\[0.50 em]
<
\!\left({f_a(c\mbox{\footnotesize $-$})
\,-\!\!\!\!\!\!
\displaystyle\sum_{\,\,\,k=\mathfrak{m}(a)+1}^{n-1}{\!\!\!E_k}{c^k}
\,-\!\!
\displaystyle\sum_{i=0}^{\mathfrak{m}(a)-1}{\!E_{i}}{c^{i}}}
\right)\!
\displaystyle\frac{{x}^n}{{c}^n}
\,\,+ \!\!\!\!
\displaystyle\sum_{\,\,\,k=\mathfrak{m}(a)+1}^{n-1}{\!\!\!E_k}{x^k}
+ \!\!
\displaystyle\sum_{i=0}^{\mathfrak{m}(a)-1}{\!E_{i}}{x^{i}},
\end{array}
\end{equation}
for every $x \!\in\! \left(0,c\right)$, $c \!\in\! \left(0,\pi\right)$, $n \!>\! \mathfrak{m}(a)+1$ and $a \in \left(\mbox{\footnotesize $\displaystyle\frac{3}{2}$},2\right)$.

\medskip
\medskip
 The family of functions $f_a(x)$,  for $x \in (0, \pi)$ and  $a \!\in\! \left(\mbox{\footnotesize $\displaystyle\frac{3}{2}$},2\right)$, satisfy the conditions $1)$ and $2)$ of Theorem~\ref{exactly_one_zero},
as we prove in the following Lemma:

%
%

\begin{lemma}\label{derivations_f}
Consider the family of functions   $
f_a(x)
=
a \ln \displaystyle\frac{\sin x}{x} - 2 \ln \cos \displaystyle\frac{x}{2}
$
for  $x \!\in\! (0,\pi)$ and  parameter  $a \!\in\! \left(\mbox{\small $\displaystyle\frac{3}{2}$},2\right)$.
Let $m=\mathfrak{m}(a)$, where $\mathfrak{m}(a)$ is defined as in~$(\ref{definition_m(a)})$.

Then, it is true that  $~\mbox{\small $\displaystyle\frac{d^k}{dx^k}$}f_{a}(x) > 0$   for $k = m, m+1, \ldots\,$
 and  $x \!\in\! (0,\pi)$,  and the following assertions hold true$:$

\medskip
\noindent
$1)$  There is a right neighbourhood  of zero in which  the following inequalities hold true$:$

 $~f_{a}(x)<0, \,\, \mbox{\small $\displaystyle\frac{d}{dx}$}f_{a}(x)<0,
\, \ldots, \mbox{\small $\displaystyle\frac{d^{m-1}}{dx^{m-1}}$}f_{a}(x)<0$,

\medskip
\noindent
$2)$ $~f_{a}(\pi_-)>0,  \,\, \mbox{\small $\displaystyle\frac{d}{dx}$}f_{a}(\pi_-)>0, \ldots,
\mbox{\small $\displaystyle\frac{d^{m-1}}{dx^{m-1}}$}f_{a}(\pi_-)>0$.
\end{lemma}

\noindent {\bf Proof.}
 Let us recall that for any fixed $a \in \left(\mbox{\footnotesize $\displaystyle\frac{3}{2}$},2\right) $ there is  a finite number of negative  coefficients in the power series expansion~(\ref{power-series-f}).
Also, we have:
$$
\left(\!\mbox{\small $\displaystyle\frac{d}{dx}$}f_{a}\!\right)\!(x)
=
a \!\left(\!\cot x \!-\! \frac{1}{x}\!\right)\! + \tan\frac{x}{2}.
$$
 For the  derivations of the function $f_{a}(x)$ in the left neighborhood  of $\pi$, it is enough to observe the following:
$$
\left(\!\mbox{\small $\displaystyle\frac{d}{dx}$}f_{a}\!\right)\!(\pi-x)
=
a \!\left(\!-\cot x \!-\! \frac{1}{\pi-x}\right) + \cot\frac{x}{2}
=
\frac{2-a}{x}-\frac{a}{\pi}+\left(\!a\!\left(\!\frac{1}{3}\!-\!\frac{1}{\pi^2}\!\right)\!-\!\frac{1}{6}\!\right)\!x+\ldots\,.
$$
From this, the conclusions of the lemma can be directly derived.  \hfill$\Box$

\bigskip
Thus, for every   $a \!\in\! \left(\mbox{\footnotesize $\displaystyle\frac{3}{2}$},2\right)$,
the corresponding function
$
f_a(x) =
a \ln \displaystyle\frac{\sin x}{x} - 2 \ln \cos \displaystyle\frac{x}{2}
$ has exactly one zero on the interval $(0, \pi)$. Let us denote it by $x_a$.

\smallskip


The following Theorem is a direct consequence of these considerations.

\begin{theorem}
\label{Corollary 11}
For every $a \!\in\! \left(\mbox{\footnotesize $\displaystyle\frac{3}{2}$},2\right)$, and every $x \in \left(0, x_a\right],$ where $0< x_a < \pi$,
the following inequality  holds true:
\begin{equation}
\label{nejednakost_Corollary11}
\left(\displaystyle\frac{\sin{x}}{x}\right)^{\!a} \leq\, \cos^{2}{\!\displaystyle\frac{x}{2}}.
\end{equation}
\end{theorem}
\smallskip

For the selected discrete values  of $a \!\in\! \left(\mbox{\footnotesize $\displaystyle\frac{3}{2}$}, 2\right) $, the zeros  $x_a$   of the corresponding functions  $f_ {a} (x)$  are  shown in  Table 1.
Although the values  $x_a$  can be obtained by any numerical method,  the following remark can also  be used to locate them.
\begin{remark}
For a fixed $a \!\in\! \left(\mbox{\footnotesize $\displaystyle\frac{3}{2}$},2\right)$, select $n>\mathfrak{m}(a)+1$
and consider inequalities~$(\ref{natural-extension-estimation})$. Denote the corresponding polynomials on the left-hand
side and the right-hand side~of $(\ref{natural-extension-estimation})$ by $P_L(x)$ and $P_R(x)$, respectively.
These polynomials are of negative sign in a right neighborhood of zero $($see$\;${\rm \cite{Malesevic2016}},$\;$Theorem$\;${\rm 2.5.}$)$,
and they have positive leading coefficients.$\,$ Then, the root $x_a$ of the equation $f_{a}(x) \!=\! 0$ is always localized
between the smallest positive root of the equation $P_L (x) \!=\! 0$ and the smallest positive root of the~equation $P_R (x) \!=\! 0$.
\end{remark}

\medskip


\subsection{Inequalities with the polynomial exponents}

In this subsection we propose and prove a new double-sided inequality involving the sinc function with polynomial exponents.

To be more specific, we find two polynomials of the second degree which, when placed in the exponent of the sinc function, give an upper and a lower bound for $ {\cos}^2 \frac{x}{2}$.

\begin{theorem}
For every $x\in (0,3.1)$ the following double-sided inequality holds$:$
\label{Teorema_PolyExp}
\begin{equation}
\label{nejednakost_PolyExp}
{\left( {\frac{{\sin x}}{x}} \right)^{\!{p_1}\left( x \right)}}
\!< \;
{\cos}^2 \frac{x}{2}
\;<\,
{\left( {\frac{{\sin x}}{x}} \right)^{\!{p_2}\left( x \right)}},
\end{equation}
where
${p_1}\left( x \right)
=
\mbox{\footnotesize $\displaystyle\frac{3}{2}$}
+
\mbox{\footnotesize $\displaystyle\frac{{{x^2}}}{{2{\pi^2}}}$}$
and
${p_2}\left( x \right)
=
\mbox{\footnotesize $\displaystyle\frac{3}{2}$}
+
\mbox{\footnotesize $\displaystyle\frac{{{x^2}}}{{80}}$}$.
\end{theorem}

{\noindent \bf Proof.} Consider the equivalent form of the inequality (\ref{nejednakost_PolyExp}):
\[{p_1}\left( x \right)\ln \frac{{\sin x}}{x} < 2\ln \cos \frac{x}{2} < {p_2}\left( x \right)\ln \frac{{\sin x}}{x}.\]

Now, let us introduce the following notation:

\[{G_i}\left( x \right) = {p_i}\left( x \right)\ln \frac{{\sin x}}{x} - 2\ln \cos \frac{x}{2},\]
for $i=1,2$.

\medskip
Based on the Theorem WD, from (\ref{Series_ln_sin_x_over_x}) we obtain:
\begin{equation}
\!\!\!\!\!\!\!\!
\begin{array}{l}
\label{lnsinxx}
-\displaystyle\sum\limits_{k=1}^{m-1}{
\mbox{\small $\displaystyle\frac{2^{2k-1}|B_{2k}|}{k(2k)!}$}x^{2k}}
+
{\Big (}\mbox{\small $\displaystyle\frac{1}{c}$}{\Big )}^{\!\!2m}\!\!
\left(
\ln \mbox{\small $\displaystyle\frac{\sin c}{c}$}
-
\displaystyle\sum\limits_{k=1}^{m-1}{
\mbox{\small $\displaystyle\frac{2^{2k-1}|B_{2k}|}{k(2k)!}$}
c^{2k}}\right)
\!x^{2m}< \\[3.0 ex]

\;\;\;\;\,
\,<\,
\ln
\mbox{\small $\displaystyle\frac{\sin x}{x}$} \,
<-\displaystyle\sum\limits_{k=1}^{n}{
\mbox{\small $\displaystyle\frac{2^{2k-1}|B_{2k}|}{k(2k)!}$}x^{2k}},
\end{array}
\end{equation}
for $x \!\in\! \left(0,\pi \right)$ where $n, m \!\in\! \NN$, $m,n \ge 2$.

\medskip

Based on the Theorem WD, from  (\ref{Series_ln_cos_x}) we obtain: 
\begin{equation}
\!\!\!\!\!\!\!\!\!\!\!\!\!\!
\begin{array}{l}
-\!\!\displaystyle\sum\limits_{k=1}^{m-1}{\!
\mbox{\small $\displaystyle\frac{2^{2k-1}(2^{2k}\!-\!1)|B_{2k}|}{k(2k)!}$}x^{2k}}
+
{\Big (}\mbox{\small $\displaystyle\frac{1}{c}$}{\Big )}^{\!2m}\!\!
\left(\!
\ln \cos c
-
\!\displaystyle\sum\limits_{k=1}^{m-1}{
\mbox{\small $\displaystyle\frac{2^{2k-1}(2^{2k}\!-\!1)|B_{2k}|}{k(2k)!}$}
c^{2k}}\!\right)
\!x^{2m}<\\[3.0 ex]
\,\,<\,\ln \cos x \, <
-\!\!\displaystyle\sum\limits_{k=1}^{n}{\!
\mbox{\small $\displaystyle\frac{2^{2k-1}(2^{2k}\!-\!1)|B_{2k}|}{k(2k)!}$}x^{2k}},
\end{array}
\!\!\!\!\!\!\!\!\!\!\!\!\!\!
\!\!\!\!\!\!\!\!\!\!\!\!\!\!
\end{equation}
for $x \!\in\! \left(0,c\right)$ and where $0 < c < \mbox{\small $\displaystyle\frac{\pi}{2}$}$, $n, m \!\in\! \NN$, $m,n \ge 2$,
i.e.
\begin{equation}
\label{lncosx2}
\!\!\!\!\!\!\!\!
\begin{array}{l}
\!\!\displaystyle\sum\limits_{k=1}^{n}{\!
\mbox{\small $\displaystyle\frac{(2^{2k}\!-\!1)|B_{2k}|}{2k(2k)!}$}x^{2k}}
<
-\ln
\cos \frac{x}{2} \, <                                                                                                  \\[3.0 ex]
\!<\!
\!\displaystyle\sum\limits_{k=1}^{m-1}{\!
\mbox{\small $\displaystyle\frac{(2^{2k}\!-\!1)|B_{2k}|}{2k(2k)!}$}x^{2k}}
-
{\Big (}\mbox{\small $\displaystyle\frac{2}{c}$}{\Big )}^{\!2m}\!\!
\left(\!
\ln \cos \frac{c}{2}
-
\!\displaystyle\sum\limits_{k=1}^{m-1}{
\mbox{\small $\displaystyle\frac{(2^{2k}\!-\!1)|B_{2k}|}{2k(2k)!}$}
c^{2k}}\!\right)
\!x^{2m},
\end{array}
\!\!\!\!\!\!\!\!\!\!\!\!\!\!
\!\!\!\!\!\!\!\!\!\!\!\!\!\!
\end{equation}
for $x \!\in\! \left(0,c\right)$ and $0 \!<\! c \!<\! \pi$, $n, m \!\in\! \NN$, $m,n \ge 2$.

\medskip
Now, let us introduce the notation:
\[\begin{array}{l}
\!\!\!\!\!\!\!\!\!\!\!\!\!\!{H_1}\left( {x,{m_1},{n_1},c_1} \right)  =  - {p_1}\left( x \right)\mathop \sum \limits_{k = 1}^{{m_1} - 1} \frac{{{2^{2k - 1}}\left| {{B_{2k}}} \right|}}{{k(2k)!}}{x^{2k}}-\\[3.0 ex]
 - 2\left( { - \mathop \sum \limits_{k = 1}^{{m_1} - 1} \frac{{\left( {{2^{2k}} - 1} \right)\left| {{B_{2k}}} \right|}}{{2k(2k)!}}{x^{2k}} + \frac{1}{{{c_1^{2m_1}}}}\left( {\ln \frac{c_1}{2} + \mathop \sum \limits_{k = 1}^{{n_1} - 1} \frac{{\left( {{2^{2k}} - 1} \right)\left| {{B_{2k}}} \right|}}{{2k(2k)!}}{c_1^{2k}}} \right){x^{2m_1}}} \right),
\end{array}\]
for $m_1,n_1\in \NN$, $m_1,n_1 \ge 2$, $c_1\in(0,\pi)$, and $x \in (0,c_1)$.

\[\begin{array}{c}
\!\!\!\!\!\!\!\!\!\!\!\!\!\!{H_2}\left( {x,{m_2},{n_2},c_2} \right) =  {p_2}\left( x \right)\left(\!\! { -\!\!\!\! \mathop \sum \limits_{k = 1}^{{m_2} - 1} \frac{{{2^{2k - 1}}\left| {{B_{2k}}} \right|}}{{k(2k)!}}{x^{2k}} + \frac{1}{{{c_2^{2m_2}}}}\left( {\ln \frac{{\sin c_2}}{c_2} + \!\!\!\! \mathop \sum \limits_{k = 1}^{{m_2} - 1} \frac{{{2^{2k - 1}}\left| {{B_{2k}}} \right|}}{{k(2k)!}}{c_2^{2k}}} \right){x^{2m_2}}} \right) +\\[3.0 ex]
+ 2\mathop \sum \limits_{k = 1}^{{n_2}} \frac{{\left( {{2^{2k}} - 1} \right)\left| {{B_{2k}}} \right|}}{{2k(2k)!}}{x^{2k}},
\end{array}\]
for $m_2,n_2\in \NN$, $m_2,n_2 \ge 2$, $c_2\in(0,\pi)$, and $x \in (0,c_2)$.\\

\break

Based on the inequalities (\ref{lnsinxx}) and (\ref{lncosx2})  the following holds true:
\[{G_1}\left( x \right) < {H_1}\left( x, m_1, n_1, c_1 \right),\]
\[{G_2}\left( x \right) > {H_2}\left( x,  m_2, n_2, c_2 \right),\]

\noindent
for $m_1,n_1,m_2,n_2 \in \NN$ and $c_1,c_2 \in (0,\pi)$.

\smallskip
For  $c_1=c_2=3.1$, $m_1=25$ and $n_1=10$ and for  $m_2=13$ and $n_2=27$,  it is easy to prove that ${H_1}\left( x, m_1, n_1, c_1 \right)<0$ and  ${H_2}\left( x, m_2, n_2, c_2 \right)>0$, for every $x \in (0,c_1)$.

\smallskip
Hence, we conclude that $G_1(x)<0$ and $G_2(x)>0$ for every $x \in (0,3.1)$, and the double-sided inequality (\ref{nejednakost_PolyExp}) holds.$ \hfill\Box$
\begin{remark}
Note that this method can be used to prove that the inequality $(\ref{nejednakost_PolyExp})$ of
Theorem \mbox{\rm \ref{Teorema_PolyExp}} holds on any interval $(0,c)$ where $c \in (0, \pi)$,
but the degrees of the polynomials $H_1$ and $H_2$ get larger as $c$ approaches $\pi$.
\end{remark}

\smallskip

\subsection{Constant vs. polynomial exponents}

Let us observe the inequalities in Theorem \ref{Corollary 11} and Theorem \ref{Teorema_PolyExp}, inequality (\ref{nejednakost_PolyExp}), containing constants and polynomials in the exponents, respectively.

\smallskip
A question of establishing a relation between these functions, with different types of exponents, comes up naturally. The following theorem addresses this question.
\begin{theorem}
\label{Teorema_const_vs_poly}
For every $a \!\in\! \left(\mbox{\footnotesize $\displaystyle\frac{3}{2}$},2\right)$ and every $x\in \left(0,m_a\right)$, where $m_a=\sqrt{2\pi^2 \left( {a - \mbox{\footnotesize $\displaystyle\frac{3}{2}$}} \right)}$, the following double-sided inequality holds$:$
\begin{equation}
\label{nejednakost_const_vs_poly}
{\left( {\frac{{\sin x}}{x}} \right)^{\!a}}
<
{\left( {\frac{{\sin x}}{x}} \right)^{\!{\frac{3}{2} + \frac{{{x^2}}}{{2{\pi ^2}}}} }} < {\cos ^2}\frac{x}{2}.
\end{equation}
\end{theorem}

\noindent{\bf Proof.}
Let $a = \mbox{\footnotesize $\displaystyle\frac{3}{2}$} + \varepsilon$,
$\varepsilon \!\in\! \left(0,\mbox{\footnotesize $\displaystyle\frac{1}{2}$}\right)$ and $x>0$. Then:
\[\begin{array}{c}
\left( {\frac{3}{2} + \frac{{{x^2}}}{{2{\pi ^2}}}} \right)\ln \frac{{\sin x}}{x} - a\ln \frac{{\sin x}}{x} = \left( {\frac{3}{2} + \frac{{{x^2}}}{{2{\pi ^2}}}} \right)\ln \frac{{\sin x}}{x} - \left( {\frac{3}{2} + \varepsilon } \right)\ln \frac{{\sin x}}{x} = \\
\\
 = \left( {\frac{{{x^2}}}{{2{\pi ^2}}} - \varepsilon } \right)\ln \frac{{\sin x}}{x} = \frac{1}{{2{\pi ^2}}}\left( {x - \sqrt {2{\pi ^2}\varepsilon } } \right)\left( {x + \sqrt {2{\pi ^2}\varepsilon } } \right)\ln \frac{{\sin x}}{x}.
\end{array}\]

Now we have:
\[\!\!\!\!\!\!\!\!
x \!\in\! \left( {0,\sqrt {2{\pi ^2}\varepsilon } } \right)
\Longleftrightarrow
\left({\mbox{\footnotesize $\displaystyle\frac{3}{2}$} + \alpha {x^2}} \right)\ln \frac{{\sin x}}{x}
>
\left( {\mbox{\footnotesize $\displaystyle\frac{3}{2}$} + \varepsilon } \right)\ln \frac{{\sin x}}{x}
\Longleftrightarrow
{\left( {\frac{{\sin x}}{x}} \right)^{\!\frac{3}{2} + \frac{{{x^2}}}{{2{\pi ^2}}}}}
\!>
{\left( {\frac{{\sin x}}{x}} \right)^{\!\frac{3}{2} + \varepsilon }}.\]

Hence, applying Theorem \ref{Teorema_PolyExp}, the double-sided inequality (\ref{nejednakost_const_vs_poly}) holds for every \break ${a \!\in\! \left(\mbox{\footnotesize $\displaystyle\frac{3}{2}$},2\right)}$ and every $x \in (0,m_a)$.
$\hfill \Box$

\break

Now, in Table 1, we show the values
  $x_a$ and $m_a$
for some specified $a \!\in\! \left( \mbox{\footnotesize $\displaystyle\frac{3}{2}$},2 \right)\,$:

{\small
$$
\!\!\!
\begin{array}{|c||c|c|c|c|c|c|c|c|c|c|c|} \hline
a   \!\!&\!\! {1.501}        \!\!&\!\! {1.502}        \!\!&\!\! {1.503}        \!\!&\!\! {1.504}        \!\!&\!\! {1.505}        \!\!&\!\! {1.506}        \!\!&\!\! {1.507}        \!\!&\!\! {1.508}        \!\!&\!\! {1.509}        \!\!&\!\! {1.510}        \\ \hline
x_a \!\!&\!\! {0.282 ...} \!\!&\!\! {0.398 ...} \!\!&\!\! {0.487 ...} \!\!&\!\! {0.561 ...} \!\!&\!\! {0.626 ...} \!\!&\!\! {0.685 ...} \!\!&\!\! {0.738 ...} \!\!&\!\! {0.788 ...} \!\!&\!\! {0.834 ...} \!\!&\!\! {0.878 ...} \\ \hline
m_a \!\!&\!\! {0.140 ...} \!\!&\!\! {0.198 ...} \!\!&\!\! {0.243 ...} \!\!&\!\! {0.280 ...} \!\!&\!\! {0.314 ...} \!\!&\!\! {0.344 ...} \!\!&\!\! {0.371 ...} \!\!&\!\! {0.397 ...} \!\!&\!\! {0.421 ...} \!\!&\!\! {0.444 ...} \\ \hline
\end{array}
$$}

\vspace*{-6.5 mm}

{\small
$$
\!\!\!
\begin{array}{|c||c|c|c|c|c|c|c|c|c|c|c|} \hline
a   \!\!&\!\! {1.52}      \!\!&\!\! {1.53}      \!\!&\!\! {1.54}      \!\!&\!\! {1.55}      \!\!&\!\! {1.56}      \!\!&\!\! {1.57}     \!\!&\!\! {1.58}       \!\!&\!\! {1.59}      \!\!&\!\! {1.60}      \!\!&\!\! {1.65}      \\ \hline
x_a \!\!&\!\! {1.220 ...} \!\!&\!\! {1.468 ...} \!\!&\!\! {1.666 ...} \!\!&\!\! {1.831 ...} \!\!&\!\! {1.973 ...} \!\!&\!\! {2.096 ...} \!\!&\!\! {2.205 ...} \!\!&\!\! {2.302 ...} \!\!&\!\! {2.302 ...} \!\!&\!\! {2.302 ...} \\ \hline
m_a \!\!&\!\! {0.628 ...} \!\!&\!\! {0.769 ...} \!\!&\!\! {0.888 ...} \!\!&\!\! {0.993 ...} \!\!&\!\! {1.088 ...} \!\!&\!\! {1.175 ...} \!\!&\!\! {1.256 ...} \!\!&\!\! {1.256 ...} \!\!&\!\! {1.256 ...} \!\!&\!\! {1.256 ...} \\ \hline
\end{array}
$$}

\vspace*{-6.5 mm}

{\small
$$
\!\!\!
\begin{array}{|c||c|c|c|c|c|c|c|c|c|c|c|} \hline
a   \!\!&\!\! {1.70}      \!\!&\!\! {1.75}      \!\!&\!\! {1.80}      \!\!&\!\! {1.85}      \!\!&\!\! {1.90}      \!\!&\!\! {1.92}      \!\!&\!\! {1.94}      \!\!&\!\! {1.96}      \!\!&\!\! {1.98}      \!\!&\!\! {1.9999}    \\ \hline
x_a \!\!&\!\! {2.911 ...} \!\!&\!\! {3.034 ...} \!\!&\!\! {3.103 ...} \!\!&\!\! {3.133 ...} \!\!&\!\! {3.141 ...} \!\!&\!\! {3.141 ...} \!\!&\!\! {3.141 ...} \!\!&\!\! {3.141 ...} \!\!&\!\! {3.141 ...} \!\!&\!\! {3.141 ...} \\ \hline
m_a \!\!&\!\! {1.986 ...} \!\!&\!\! {2.221 ...} \!\!&\!\! {2.433 ...} \!\!&\!\! {2.628 ...} \!\!&\!\! {2.809 ...} \!\!&\!\! {2.879 ...} \!\!&\!\! {2.947 ...} \!\!&\!\! {3.013 ...} \!\!&\!\! {3.087 ...} \!\!&\!\! {3.141 ...} \\ \hline
\end{array}
$$}

\vspace*{-5.0 mm}

$$ \mbox{\bf Table 1} $$

\begin{remark}
Note that Theorem $\ref{Teorema_const_vs_poly}$ represents another proof of the following assertion from $\mbox{\rm \cite{Lutovac2017}}\!:$
$$
{\big (}\,\forall a \!\in\! \left(3/2,2\right)\!{\big )}
{\big (}\,\exists \delta \!>\! 0\,{\big )}
{\big (}\,\forall x \!\in\! (0, \delta){\big )}
{\left( {\frac{{\sin x}}{x}} \right)^{\!\!a}} \!\!< {\cos ^2}\frac{x}{2}.
$$
\end{remark}
\medskip

\section{Conclusion}

In this paper,  using the power series expansions and the application of the {\sc Wu}-{\sc Debnath} theorem, we proved that the inequality (\ref{Z-H-Jang})
holds for $a \!=\! \mbox{\footnotesize $\displaystyle\frac{3}{2}$}$. At the same time, this proof represents  a new short proof of  Statement~\ref{Brankova-teorema}.

We analyzed the cases $a \!\in\! \left(\mbox{\footnotesize $\displaystyle\frac{3}{2}$},2\right)$ and $\displaystyle a \!\geq\! 2$
and we prove  the corresponding inequalities.
We introduced and prove a new double-sided inequality of similar type involving polynomial
exponents.  Also, we established a relation between the cases of the constant  and of the polynomial  exponent.

\bigskip

\bigskip
\noindent \textbf{Acknowledgement.} The research of the first, second and third authors was supported
in part by the Serbian Ministry of Education, Science and Technological Development, under Projects
ON 174033, TR 32023, and ON 174032 \& III 44006, respectively.

\bigskip
\noindent \textbf{Competing Interests.} The authors would like to
state that they do not have any competing interests in the subject
of this research.

\bigskip
\noindent \textbf{Author's Contributions.} All the authors
participated in every phase of the research conducted for this
paper.

\break

\end{document}